\documentclass[final]{siamltex}

\usepackage{amsmath}
\usepackage{amssymb}

\title{Convergence acceleration algorithm via an equation related to the
  lattice Boussinesq equation}

\author{Yi He\footnotemark[2] \footnotemark[3] , Xing-Biao
  Hu\footnotemark[2] , Jian-Qing Sun\footnotemark[2] \footnotemark[3] ,
  \and Ernst Joachim Weniger\footnotemark[4]}

\begin{document}
\maketitle

\renewcommand{\thefootnote}{\fnsymbol{footnote}}
\footnotetext[2]{LSEC, Institute of Computational Mathematics and
  Scientific Engineering Computing, AMSS, Chinese Academy of Sciences,
  P.O.Box 2719, Beijing 100190, PR China.}
\footnotetext[3]{Graduate School of the Chinese Academy of
  Sciences,\ Beijing, PR China.}
\footnotetext[4]{Institut f\"ur Physikalische und Theoretische
  Chemie, Universit\"at Regensburg, D-8400 Regensburg, Federal Republic
  of Germany.}

\renewcommand{\thefootnote}{\arabic{footnote}}



\begin{abstract}
  The molecule solution of an equation related to the lattice Boussinesq
  equation is derived with the help of determinantal identities. It is
  shown that this equation can for certain sequences be used as a
  numerical convergence acceleration algorithm. Numerical examples with
  applications of this algorithm are presented.
\end{abstract}

\begin{keywords}
  Lattice Boussinesq equation, convergence acceleration algorithm,
  molecule solution
\end{keywords}

\begin{AMS}
  65B05, 37K40, 37K60
\end{AMS}

\pagestyle{myheadings}
\thispagestyle{plain}

\section{Introduction}
\label{Sec:Intro}

In recent years, it has been found that integrable systems are closely
connected to certain numerical algorithms. This observation allows a
fresh look at the research in both fields, and a lot of interesting work
has been done lately.

For example, one step of the QR algorithm is equivalent to the time
evolution of the finite nonperiodic Toda lattice \cite{16}.  Moreover,
Wynn's celebrated $\varepsilon$-algorithm \cite{19} is nothing but the
fully discrete potential KdV equation \cite{10, 14}.  The continuous-time
Toda equation leads to a new algorithm for computing the Laplace
transform of a given analytic function \cite{11}. The discrete
Lotka--Volterra system has applications in numerical algorithms for
computing singular values \cite{7, 8, 18}.  As far as the links between
integrable systems and convergence acceleration algorithms are concerned,
more results have recently been achieved (see, e.g., \cite{3a, 5a, 9a,
  9b}).

On the other hand, rapid progress has been made recently in the study of
discrete integrable systems. As a result, many new examples of discrete
integrable systems have been found, among them the lattice
Gel'fand--Dikii hierarchy \cite{12, 13}. However, to the best of our
knowledge, nothing has been done so far on designing convergence
acceleration algorithms via new discrete integrable systems. This is what
we want to do in this article. Our starting point will be the lattice
Boussinesq equation \cite{12, 13}, and we want to use it for the design
of a new convergence acceleration algorithm.

In view of the importance of both discrete integrable systems and
sequence transformations in the context of convergence acceleration
algorithms, it might be helpful -- before discussing more specific
details -- to first mention some basic facts about discrete integrable
systems and about sequence transformations, and to explain how discrete
integrable systems can be a source of inspiration for sequence
transformations.

Discrete integrable systems could be considered as a specific class of
discrete systems which possesses an important property of what is usually
called \emph{integrability}. In particular, the existence of a Lax pair
and a $\tau$-function are two of the most important features of
integrability shared by some famous numerical algorithms, such as
Rutishauser's $qd$-algorithm \cite{14b} and Wynn's
$\varepsilon$-algorithm \cite{19}. Let us also remark that corresponding
to different boundary conditions there are different versions of discrete
integrable systems available. An example is the famous Toda equation
\begin{equation}
  \label{TodaLatEq}
  \frac {d^2x_k}{dt^2}=e^{x_{k-1}-x_k}-e^{x_{k}-x_{k+1}}.
\end{equation}
In the infinite chain case with $k=0,\pm 1,\pm 2, \ldots $, we call
(\ref{TodaLatEq}) the \emph{infinite Toda lattice equation}. Under
periodic boundary condition $x_{k+K_0}=x_k$ with $k=0,\pm 1,\pm 2,
\ldots$ and fixed $K_{0} > 0$, (\ref{TodaLatEq}) is referred to as the
\emph{periodic Toda lattice equation}. If $k=0, 1, 2, \ldots$ with the
boundary condition $x_0(t)=-\infty$, we call (\ref{TodaLatEq}) the
\emph{semi-infinite Toda equation} or \emph{infinite Toda molecule
  equation}. If $x_0(t)=-\infty$ and $x_{N+1}=+\infty $, we call
(\ref{TodaLatEq}) the \emph{finite nonperiodic Toda equation} or
\emph{finite molecule Toda equation}.

In connection with convergence acceleration algorithms, we are
exclusively interested in the semi-infinite or infinite molecule case of
(\ref{TodaLatEq}), corresponding to the semi-infinite or infinite Toda
molecule equation, and the solutions obtained in this way are called
molecule solutions. As it turns out, our molecule solutions are closely
related to sequence transformations.

We now return to the lattice Boussinesq equation, which is the second
equation of the lattice Gel'fand--Dikii hierarchy \cite{12, 13}. The
multisoliton solutions of the lattice Boussinesq equation defined on the
elementary square were derived in \cite{6, 17}. An ultradiscrete lattice
Boussinesq equation and an alternate form of the discrete potential
Boussinesq equation had been proposed and the multisoliton solutions of
both equations had also been obtained \cite{9}. The lattice Boussinesq
equation is expressed as \cite{12, 13}
\begin{align*}
  & \frac{p^3-q^3}{p-q+u_{l+1}^{m+1}-u_{l+2}^{m}} -
  \frac{p^3-q^3}{p-q+u_{l}^{m+2}-u_{l+1}^{m+1}}
  \\
  & \quad - u_{l}^{m+1} u_{l+1}^{m+2} + u_{l+1}^{m} u_{l+2}^{m+1} +
  u_{l+2}^{m+2} (p-q+u_{l+1}^{m+2} - u_{l+2}^{m+1})
  \\
  & \quad \quad + u_{l}^{m} (p-q+u_{l}^{m+1}-u_{l+1}^{m})
  \\
  & \quad \quad \quad \; = \; (2p+q) (u_{l+1}^{m} + u_{l+1}^{m+2}) -
  (p+2q) (u_{l}^{m+1}+u_{l+2}^{m+1}) ,
\end{align*}
which is equivalent to
\begin{align}
  \label{bsq1}
  & \frac{p^3-q^3}{p-q+u_{l+1}^{m+1} - u_{l+2}^{m}} -
  \frac{p^3-q^3}{p-q+u_{l}^{m+2}-u_{l+1}^{m+1}}
  \notag \\
  & \quad \; = \; (p+2q+u_{l+1}^{m}-u_{l+2}^{m+2})
  (p-q+u_{l+1}^{m+2}-u_{l+2}^{m+1})
  \notag \\
  & \qquad \quad - (p+2q+u_{l}^{m}-u_{l+1}^{m+2})
  (p-q+u_{l}^{m+1}-u_{l+1}^{m}) .
\end{align}
If we now set
\[
\bar{u}_l^m=u_l^m-pl-qm, \qquad p^3-q^3=1,
\]
we obtain from \eqref{bsq1}
\begin{align}
  \label{bsq2}
  & \frac{1}{\bar{u}_{l+1}^{m+1}-\bar{u}_{l+2}^{m}} -
  \frac{1}{\bar{u}_{l}^{m+2}-\bar{u}_{l+1}^{m+1}}
  \notag \\
  & \quad \; = \; (\bar{u}_{l+1}^{m}-\bar{u}_{l+2}^{m+2})
  (\bar{u}_{l+1}^{m+2}-\bar{u}_{l+2}^{m+1}) -
  (\bar{u}_{l}^{m}-\bar{u}_{l+1}^{m+2})
  (\bar{u}_{l}^{m+1}-\bar{u}_{l+1}^{m}) .
\end{align}
With the help of the variable transformations
\[
n=-l, \qquad k=m+l, \qquad U_k^n=\bar{u}_l^m,
\]
we obtain from \eqref{bsq2}
\begin{align*}
  & (U_{k+1}^{n+1}-U_{k+4}^{n}) (U_{k+3}^{n+1}-U_{k+3}^{n}) -
  (U_k^{n+2}-U_{k+3}^{n+1})(U_{k+1}^{n+2}-U_{k+1}^{n+1})
  \\
  & \quad \; = \; \frac{1}{U_{k+2}^{n+1}-U_{k+2}^{n}} -
  \frac{1}{U_{k+2}^{n+2}-U_{k+2}^{n+1}} .
\end{align*}
This relationship can be simplified further yielding the following
equation:
\begin{equation}
  \label{ag1}
  U_{k+3}^n = U_k^{n+1} - \frac{1}{(U_{k+2}^{n+1} -
    U_{k+2}^{n})(U_{k+1}^{n+1}-U_{k+1}^{n})} .
\end{equation}

In this article, we will first derive the molecule solution of the
two-dimensional difference equation \eqref{ag1}. Then we will show that
the resulting equation can be used for the acceleration of the
convergence of computationally relevant sequences.

Our article is organized as follows: In section \ref{Sec:MolSolBoussEq},
we will derive the molecule solution of \eqref{ag1} with the help of
determinantal identities. In section \ref{Sec:GenPropSeqTr}, we will
provide a highly condensed review of the most basic features of sequence
transformations. In section \ref{Sec:ConvAccAlg}, we will use the results
from section \ref{Sec:MolSolBoussEq} to construct a new sequence
transformation. We will also show that this transformation can be
implemented by the lattice equation \eqref{ag1} with given initial
values. In section \ref{Sec:NumExp}, applications of this algorithm are
presented. Section \ref{Sec:ConclDisc} is devoted to conclusions and
discussions.

\section{Molecule solution of the lattice equation \eqref{ag1}}
\label{Sec:MolSolBoussEq}

In this section, we study the molecule solution of the lattice
equation \eqref{ag1} by Hirota's bilinear method
\cite{7a}. With the help of the dependent variable transformation
\[
U_k^{n}=\frac{G_{k}^n}{F_{k}^n},
\]
we obtain the bilinear form of \eqref{ag1}
\begin{align}
  \label{be1}
  F_k^nG_k^{n+1} - F_k^{n+1} G_k^{n} & \; = \; F_{k+1}^n F_{k-1}^{n+1} ,
  \\
  \label{be2}
  F_{k+2}^n G_{k-1}^{n+1} - F_{k-1}^{n+1}G_{k+2}^{n} & \; = \;
  F_{k+1}^{n+1} F_{k}^{n} ,
  \\
  \label{be3}
  F_{k}^nF_{k+2}^{n+1} - F_{k}^{n+1} F_{k+2}^{n} & \; = \; F_{k+3}^n
  F_{k-1}^{n+1} .
\end{align}
Equation \eqref{be3} can be derived from \eqref{be1}--\eqref{be2} by
eliminating the $G$'s.

Set
\begin{align*}
  \Psi_k (v_n) & \; = \;
  \begin{vmatrix}
    v_n              &  v_{n+1}            & \cdots  & v_{n+k-1}\\
    \Delta^2 v_n     & \Delta^2 v_{n+1}    & \cdots  & \Delta^2 v_{n+k-1}\\
    \Delta^4 v_n     & \Delta^4 v_{n+1}    & \cdots  & \Delta^4 v_{n+k-1}\\
    \vdots & \vdots & & \vdots\\
    \Delta^{2k-2} v_n & \Delta^{2k-2} v_{n+1} & \cdots & \Delta^{2k-2}
    v_{n+k-1}
  \end{vmatrix}, \qquad k=1, 2, \ldots, \\
  \Psi_{-1}(v_n) & \; = \; 0 , \quad \Psi_0(v_n)=1 ,
\end{align*}
and
\begin{align*}
  \Phi_k(v_n)& \; = \;
  \begin{vmatrix}
    n            &  n+1                & \cdots  & n+k-1 \\
    v_n          &  v_{n+1}            & \cdots  & v_{n+k-1}\\
    \Delta^2 v_n & \Delta^2 v_{n+1}    & \cdots  & \Delta^2 v_{n+k-1}\\
    \Delta^4 v_n & \Delta^4 v_{n+1}    & \cdots  & \Delta^4 v_{n+k-1}\\
    \vdots       & \vdots              &         & \vdots\\
    \Delta^{2k-4} v_n & \Delta^{2k-4} v_{n+1} & \cdots & \Delta^{2k-4}
    v_{n+k-1}
  \end{vmatrix}, \qquad k=1, 2, \ldots, \\
  \Phi_{-1}(v_n)& \; = \; 0, \quad \Phi_0 (v_n) = 1.
\end{align*}
The solution of an initial value problem related to equations \eqref{be1}
-- \eqref{be3} is given below.\\
\\
\textbf{Theorem 1} Given the initial values
\begin{equation}
  \label{cd}
  F_1^n = F_2^n = F_3^n = 1,
  \quad F_4^n = \Delta S_n,\quad G_1^n = 0,\quad G_2^n = n,\quad G_3^n = S_n.
\end{equation}
the solution of the bilinear equations \eqref{be1}--\eqref{be3} can
be expressed as follows:
\begin{align}
  \label{s1}
  F_{3k}^n & \; = \; \Psi_{k-1} (\Delta^3 S_n), \ & G_{3k}^n & =
  \Psi_k(S_n),
  \\
  F_{3k+1}^n & \; = \; \Psi_{k} (\Delta S_n), \ & G_{3k+1}^n & = -
  \Psi_{k-1}(\Delta^4S_n),
  \\
  \label{s2}
  F_{3k+2}^n & \; = \; \Psi_{k} (\Delta^2S_n), \ & G_{3k+2}^n &
  =\Phi_{k+1}(\Delta S_n).
\end{align}
Proof: We consider the case $k=3m$ in \eqref{be1}--\eqref{be3}. First, we
prove the validity of the relationship
\begin{equation}\label{cq1}
  F_{3m}^nG_{3m}^{n+1}-F_{3m}^{n+1}G_{3m}^{n}=F_{3m+1}^nF_{3m-1}^{n+1}.
\end{equation}
Define
\[ D_1 \equiv
\begin{vmatrix}
  1               &  1                 & \cdots  & 1\\
  S_n             &  S_{n+1}           & \cdots  & S_{n+m}\\
  \Delta^2 S_n    & \Delta^2 S_{n+1}   & \cdots  & \Delta^2 S_{n+m}\\
  \vdots & \vdots & &\vdots\\
  \Delta^{2m-2} S_n & \Delta^{2m-2} S_{n+1} & \cdots & \Delta^{2m-2}
  S_{n+m}
\end{vmatrix}.
\]
Then we have the relations
\begin{eqnarray}
  &&D_1=\Psi_m(\Delta S_n)=F_{3m+1}^n,\\
  &&D_1(1, 2|1, m+1)=\Psi_{m-1}(\Delta^2S_{n+1})=F_{3m-1}^{n+1},\\
  &&D_1(1|1)=\Psi_{m}(S_{n+1})=G_{3m}^{n+1},\\
  &&D_1(2|m+1)=\Psi_{m-1}(\Delta^3S_{n})=F_{3m}^{n},\\
  &&D_1(1|m+1)=\Psi_{m}(S_{n})=G_{3m}^{n},\\
  &&D_1(2|1)=\Psi_{m-1}(\Delta^3S_{n+1})=F_{3m}^{n+1},
\end{eqnarray}
where $D_1(j|k)$ and
$D_1(j, k|p, q)$ are $m$th-order and
$(m-1)$th-order determinants obtained by eliminating the $j$th row and
the $k$th column from the $D_1$ and by eliminating the $j$th and $k$th
rows and the $p$th and $q$th columns from the determinant $D_1$,
respectively.\\
From the above results, we see that the bilinear equation
\eqref{cq1} is equivalent to the Jacobi identity \cite{4}
\[
D_1D_1(1, 2|1, m+1)= D_1(1|1)D_1(2|m+1)-D_1(1|m+1)D_1(2|1).
\]
Next, we prove the validity of another relationship
\begin{equation}\label{cq2}
  F_{3m+2}^nG_{3m-1}^{n+1}-F_{3m-1}^{n+1}G_{3m+2}^{n} =
  F_{3m+1}^{n+1}F_{3m}^{n}.
\end{equation}
According to the assumptions of theorem 1, we have
\begin{eqnarray*}
lhs&=&F_{3m+2}^n G_{3m-1}^{n+1} - F_{3m-1}^{n+1} G_{3m+2}^{n}
  \\
  &=&\begin{vmatrix}
    1               &  1                 & \cdots  & 1\\
    \Delta S_n      &\Delta S_{n+1}      & \cdots  & \Delta S_{n+m}\\
    \Delta^3 S_n     & \Delta^3 S_{n+1}    & \cdots  & \Delta^3 S_{n+m}\\
    \vdots & \vdots & &\vdots\\
    \Delta^{2m-1} S_n & \Delta^{2m-1} S_{n+1} & \cdots & \Delta^{2m-1}
    S_{n+m}
  \end{vmatrix}
  \\[2.0\jot]
  &\times&
  \begin{vmatrix}
    n+1                  & n+2                & \cdots  & n+m\\
    \Delta S_{n+1}       &\Delta S_{n+2}      & \cdots  & \Delta S_{n+m}\\
    \Delta^3 S_{n+1}      & \Delta^3 S_{n+2}    & \cdots  & \Delta^3 S_{n+m}\\
    \vdots & \vdots & &\vdots\\
    \Delta^{2m-3} S_{n+1} & \Delta^{2m-3} S_{n+2} & \cdots &
    \Delta^{2m-3} S_{n+m}
  \end{vmatrix}
  \\[2.0\jot]
  & -&
  \begin{vmatrix}
    1               &  1                 & \cdots  & 1\\
    \Delta S_{n+1}  &\Delta S_{n+2}      & \cdots  & \Delta S_{n+m}\\
    \Delta^3 S_{n+1} & \Delta^3 S_{n+2}    & \cdots  & \Delta^3 S_{n+m}\\
    \vdots & \vdots & &\vdots\\
    \Delta^{2m-3} S_{n+1} & \Delta^{2m-3} S_{n+2} & \cdots &
    \Delta^{2m-3} S_{n+m}
  \end{vmatrix}
   \\[1.0\jot]
   & \times &
  \begin{vmatrix}
    n                  & n+1                & \cdots  & n+m\\
    \Delta S_{n}       &\Delta S_{n+1}      & \cdots  & \Delta S_{n+m}\\
    \Delta^3 S_{n}      & \Delta^3 S_{n+1}    & \cdots  & \Delta^3 S_{n+m}\\
    \vdots & \vdots & &\vdots\\
    \Delta^{2m-1} S_{n} & \Delta^{2m-1} S_{n+1} & \cdots & \Delta^{2m-1}
    S_{n+m}
  \end{vmatrix}.
\end{eqnarray*}
By Schwein's determinantal identity \cite{1} we obtain
\begin{eqnarray*}
  lhs &=&\begin{vmatrix}
    1               &  1                 & \cdots  & 1\\
    n               & n+1                & \cdots  & n+m\\
    \Delta S_n      &\Delta S_{n+1}      & \cdots  & \Delta S_{n+m}\\
    \Delta^3 S_n     & \Delta^3 S_{n+1}    & \cdots  & \Delta^3 S_{n+m}\\
    \vdots & \vdots & &\vdots\\
    \Delta^{2m-3} S_n & \Delta^{2m-3} S_{n+1} & \cdots & \Delta^{2m-3}
    S_{n+m}
  \end{vmatrix}\\
  &\times&
  \begin{vmatrix}
    \Delta S_{n+1}       &\Delta S_{n+2}      & \cdots  & \Delta S_{n+m}\\
    \Delta^3 S_{n+1}      & \Delta^3 S_{n+2}    & \cdots  & \Delta^3 S_{n+m}\\
    \vdots & \vdots & &\vdots\\
    \Delta^{2m-1} S_{n+1} & \Delta^{2m-1} S_{n+2} & \cdots &
    \Delta^{2m-3} S_{n+m}
  \end{vmatrix}.
\end{eqnarray*}
We obviously have
\begin{eqnarray*}
  &&\begin{vmatrix}
    1               &  1                 & \cdots  & 1\\
    n               & n+1                & \cdots  & n+m\\
    \Delta S_n      &\Delta S_{n+1}      & \cdots  & \Delta S_{n+m}\\
    \Delta^3 S_n     & \Delta^3 S_{n+1}    & \cdots  & \Delta^3 S_{n+m}\\
    \vdots & \vdots & &\vdots\\
    \Delta^{2m-3} S_n & \Delta^{2m-3} S_{n+1} & \cdots & \Delta^{2m-3}
    S_{n+m}
  \end{vmatrix}
  \end{eqnarray*}
 \begin{eqnarray*}
  &=&\begin{vmatrix}
    1               &  1                 & \cdots  & 1\\
    \Delta^2 S_n    &\Delta^2 S_{n+1}    & \cdots  & \Delta^2 S_{n+m-1}\\
    \Delta^4 S_n     & \Delta^4 S_{n+1}    & \cdots  & \Delta^4 S_{n+m-1}\\
    \vdots & \vdots & &\vdots\\
    \Delta^{2m-2} S_n & \Delta^{2m-2} S_{n+1} & \cdots & \Delta^{2m-2}
    S_{n+m}
  \end{vmatrix}
 \\[2\jot]
&=&\begin{vmatrix}
    \Delta^3 S_n    &\Delta^3  S_{n+1}    & \cdots  & \Delta^3 S_{n+m-2}\\
    \Delta^5 S_n    & \Delta^5 S_{n+1}    & \cdots  & \Delta^5 S_{n+m-2}\\
    \vdots & \vdots & &\vdots\\
    \Delta^{2m-1} S_n & \Delta^{2m-1} S_{n+1} & \cdots & \Delta^{2m-1}
    S_{n+m-2}
  \end{vmatrix}
  \\[2\jot]
  &=& \Psi_{m-1} (\Delta^3 S_n) = F_{3m}^n.
\end{eqnarray*}
Thus, $lhs = F_{3m}^{n}F_{3m+1}^{n+1}$, and \eqref{cq2} is
proved.\\
Next, we prove the third relationship:
\begin{equation}\label{cq3}
  F_{3m}^n F_{3m+2}^{n+1} - F_{3m}^{n+1} F_{3m+2}^{n} =
  F_{3m+3}^n F_{3m-1}^{n+1}.
\end{equation}
Define
\[D_2 \equiv
\begin{vmatrix}
  1                &  1                  & \cdots  & 1\\
  \Delta^2 S_n     & \Delta^2 S_{n+1}    & \cdots  & \Delta^2 S_{n+m}\\
  \Delta^4 S_n     & \Delta^4 S_{n+1}    & \cdots  & \Delta^4 S_{n+m}\\
  \vdots & \vdots & &\vdots\\
  \Delta^{2m} S_n & \Delta^{2m} S_{n+1} & \cdots & \Delta^{2m} S_{n+m}
\end{vmatrix}.
\]
Then we have the relations
\begin{eqnarray}
  &&D_2=\Psi_m(\Delta^3S_n)=F_{3m+3}^n,\\
  &&D_2(1, m+1|1, m+1)=\Psi_{m-1}(\Delta^2S_{n+1})=F_{3m-1}^{n+1},\\
  &&D_2(1|1)=\Psi_{m}(\Delta^2S_{n+1})=F_{3m+2}^{n+1},\\
  &&D_2(m+1|m+1)=\Psi_{m-1}(\Delta^3S_{n})=F_{3m}^{n},\\
  &&D_2(1|m+1)=\Psi_{m}(\Delta^2S_{n})=F_{3m+2}^{n},\\
  &&D_2(m+1|1)=\Psi_{m-1}(\Delta^3S_{n+1})=F_{3m}^{n+1}.
\end{eqnarray}
From the above results, we see that the bilinear equation \eqref{cq3} is
nothing but the Jacobi identity \cite{4}
\[
D_2D_2(1, m+1|1, m+1)=
D_2(1|1)D_2(m+1|m+1)-D_2(1|m+1)D_2(m+1|1).
\]
The proof of other cases of \eqref{be1}--\eqref{be3} can be
obtained in a similar way. $\Box$

\section{General properties of sequence transformations}
\label{Sec:GenPropSeqTr}

Many calculations produce results that are actually sequences whose rate
of convergence is governed by one or several parameters. Unfortunately,
it often happens that the rate of convergence of such a sequence $\{
S_{n} \}$ is so slow that the determination of a sufficiently accurate
approximation to its limit $S = S_{\infty}$ by increasing the index $n$
does not work in practice. Another frequently occurring problem is that
such a sequence $\{ S_{n} \}$ does not necessarily produce a convergent
result as $n \to \infty$ even if it actually corresponds to a meaningful
mathematical object with a well-defined numerical value.

In such a situation, it can be extremely helpful to apply a so-called
\emph{sequence transformation}, which transforms the original sequence
$\{ S_n \}$ to a new sequence $\{ S'_{n} \}$ with hopefully better
numerical properties according to
\begin{displaymath}
  \mathcal{T} \colon \{ S_{n} \} \longmapsto \{ S'_{n} \} \, .
\end{displaymath}
In rudimentary form, sequence transformations have been known for
centuries. Their older history is reviewed in an article \cite{2e} and a
monograph \cite{2a} by Brezinski. More recent developments are discussed
in two articles by Brezinski \cite{2c,2d}. There is also an extensive
bibliography by Brezinski \cite{2b} containing more than 6000 references
up to 1991.

The active research on sequence transformations is documented by the fact
that in recent years quite a few specialized monographs or longer reviews
have been published, for example the ones by Brezinski and Redivo Zaglia
\cite{3}, Sidi \cite{15a}, and Weniger \cite{18a}. Numerous other
references can be found in \cite[Appendix B]{18c}.

Sequence transformations try to achieve an acceleration of convergence or
a summation in the case of divergence by purely numerical means. Since,
however, a computational algorithm can involve only a finite number of
arithmetic operations, a sequence transformation $\mathcal{T}$ can
associate only a \emph{finite subset} of the input sequence $\{ S_{n} \}$
with an element $S'_{m}$ of the transformed sequence.

All the commonly occurring sequence transformations $\mathcal{T}$ can be
represented by an infinite set of doubly indexed quantities $T_k^{(n)}$
with $k, n \in \mathbb{N}_0$. In the literature, the superscript $n$
typically indicates the minimal index occurring in the finite subset $\{
S_n, S_{n+1}, \ldots , S_{n+\ell} \}$ with $\ell = \ell (k)$ of sequence
elements which are used for the computation of the transform $T_k^{(n)}$,
and the subscript $k$ -- usually called the \emph{order} of the
transformation -- is a measure for the complexity of
$T_k^{(n)}$. Moreover, the $T_k^{(n)}$ are gauged in such a way that
$T_0^{(n)}$ always corresponds to an untransformed sequence element
according to $T_0^{(n)} = S_n$.

The basic assumption of all sequence transformations is that the elements
of a slowly convergent or divergent sequence $\{ S_n \}$, which could be
the partial sums $S_n = \sum_{k=0}^{n} a_{k}$ of an infinite series, can
for all indices $n$ be partitioned into a (generalized) limit $S$ and a
remainder or truncation error $R_n$ according to $S_{n} = S + R_n$. If
the sequence $\{ S_n \}$ converges to its limit $S$, the remainders $R_n$
can be made negligible by increasing $n$ as much as necessary. But many
sequences converge so slowly that this does not work in
practice. Increasing the index $n$ also does not help in the case of a
divergent sequence.

Alternatively, one can try to improve convergence or accomplish a
summation by computing approximations to the remainders $R_{n}$ which are
then eliminated from the sequence elements $S_{n}$, yielding a new
sequence with elements $S'_{n} = S + R'_{n}$. At least conceptually, this
is what a sequence transformation tries to do.

Some transformations -- for example the so-called Levin-type
transformations discussed in \cite{18b} -- possess closed form
expressions. However, the vast majority of all known sequence
transformations are defined by a recursive scheme.

The probably best known example of such a transformation is Wynn's
$\varepsilon$-algorithm \cite{19}, which corresponds to the following
nonlinear two-dimensional recursive scheme:
\begin{subequations}
  \label{eps_al}
  \begin{align}
    \label{eps_al_a}
    \varepsilon_{-1}^{(n)} & \; = \; 0 \, ,
    \qquad \varepsilon_0^{(n)} \, = \, S_n \, ,
    \qquad  n \in \mathbb{N}_0 \, , \\
    \label{eps_al_b}
    \varepsilon_{k+1}^{(n)} & \; = \; \varepsilon_{k-1}^{(n+1)} \, + \,
    \frac{1}{\varepsilon_{k}^{(n+1)} - \varepsilon_{k}^{(n)}} \, ,
    \qquad k, n \in \mathbb{N}_0 \, .
  \end{align}
\end{subequations}
The elements $\varepsilon_{2k}^{(n)}$ with \emph{even} subscripts provide
approximations to the (generalized) limit $S$ of the sequence $\{ S_n \}$
to be transformed, whereas the elements $\varepsilon_{2k+1}^{(n)}$ with
\emph{odd} subscripts are only auxiliary quantities which diverge if the
whole process converges. A compact FORTRAN 77 program for the
$\varepsilon$-algorithm as well as the underlying computational algorithm
is described in \cite[section 4.3]{18a}, and in \cite[p.\ 213]{14a}, a
translation of this FORTRAN 77 program to \texttt{C} can be found.

The two-dimensional recursive scheme (\ref{eps_al}) for the
$\varepsilon$-algorithm was derived by making some assumptions about the
mathematical nature of the truncation errors or remainders $R_{n} = S_{n}
- S$ of the elements of the sequence $\{ S_{n} \}$ which is to be
transformed. Other sequence transformations defined by recursive schemes
were also constructed by trying to eliminate the remainders of
appropriate model sequences.

But a converse approach is also possible. We can take an equation in two
\emph{discrete} variables -- for example, one based on integrable systems
-- and analyze if and under which conditions this equation can be
used as a starting point for the construction of a sequence
transformation.

Let us conclude this section with some useful terminology that is typical
of the literature on sequence transformations. Assume that a sequence $\{
S_n \}$, which converges to some limit $S$, satisfies
\begin{equation}
\lim_{n \to \infty} \; \frac {S_{n+1} - S} {S_n - S} \; = \;
\lim_{n \to \infty} \; \frac {R_{n+1}}{R_n} \; =\; \rho \, .
\label{DefLinLogConv}
\end{equation}
If $0 < \vert \rho \vert < 1$ holds, we say that the sequence $\{ S_n \}$
converges \emph{linearly}; if $\rho = 1$ holds, we say that this sequence
converges \emph{logarithmically}; and if $\rho = 0$ holds, we say that it
converges \emph{hyperlinearly}. Of course, $\vert \rho \vert > 1$ implies
that the sequence $\{ S_n \}$ diverges. Simple examples of linearly,
logarithmically, and hyperlinearly convergent sequences are the partial
sums of the geometric series $1/(1-z) = \sum_{\nu=0}^{\infty} z^{\nu}$,
of the Dirichlet series $\zeta (s) = \sum_{\nu=0}^{\infty} (\nu+1)^{-s}$
for the Riemann zeta function, and of the power series $\exp (z) =
\sum_{\nu=0}^{\infty} z^{\nu}/\nu!$ for the exponential function,
respectively.

Let us assume that two sequences $\{ S_n \}$ and $\{ S^{\prime}_n \}$
converge to the same limit $S$. We say that the sequence $\{ S^{\prime}_n
\}$ converges \emph{more rapidly} than $\{ S_n \}$ if
\begin{equation}
  \label{Def_AccConv}
  \lim_{n \to \infty} \; \frac { S^{\prime}_n - S }{ S_n - S}
  \; = \; 0 \, .
\end{equation}

\section{A convergence acceleration algorithm}
\label{Sec:ConvAccAlg}

In this section, we propose a new sequence transformation and show
that this transformation can be derived via the lattice
equation \eqref{ag1}.

We now consider a new sequence transformation defined as the following
ratio of determinants:
\begin{equation}
  \label{tr}
  T_k^{(n)}=\frac{\begin{vmatrix}
      S_n             &  S_{n+1}           & \cdots  & S_{n+k}\\
      \Delta^2 S_n    & \Delta^2 S_{n+1}   & \cdots  & \Delta^2 S_{n+k}\\
      \Delta^4 S_n    & \Delta^4 S_{n+1}   & \cdots  & \Delta^4 S_{n+k}\\
      \vdots & \vdots & &\vdots\\
      \Delta^{2k} S_n  & \Delta^{2k} S_{n+1} & \cdots  & \Delta^{2k} S_{n+k}
    \end{vmatrix}} {\begin{vmatrix}
      1               & 1                  & \cdots  & 1  \\
      \Delta^2 S_n    & \Delta^2S_{n+1}   & \cdots  & \Delta^2 S_{n+k}\\
      \Delta^4 S_n    & \Delta^4S_{n+1}   & \cdots  & \Delta^4 S_{n+k}\\
      \vdots & \vdots & &\vdots\\
      \Delta^{2k} S_n & \Delta^{2k}S_{n+1} & \cdots & \Delta^{2k} S_{n+k}
    \end{vmatrix}
  }.
\end{equation}
Obviously this transformation is a particular case of the
$E$-transformation which was first derived by Schneider \cite{15}
and later rederived independently by H{\aa}vie \cite{5} and
Brezinski \cite{2}. From the kernel of the $E$-transformation
\cite{3}, the kernel of the transformation \eqref{tr} is obtained in
the theorem below.

\noindent \textbf{Theorem 2} A necessary and sufficient condition that
$T_k^{(n)}=S$ for all $n$ is that
\[
S_n=S+a_1\Delta^2S_n+a_2\Delta^4S_n+\cdots+a_{k}\Delta^{2k}S_n,
\]
where $a_i$ with $ i=1,\ldots,k$ are some constants.

Theorem 1 implies that the transformation \eqref{tr} can be implemented
via the lattice equation
\begin{equation}\label{am1}
  U_{k+3}^n = U_k^{n+1} - \frac{1}{(U_{k+2}^{n+1} -
    U_{k+2}^{n})(U_{k+1}^{n+1} - U_{k+1}^{n})},
  \qquad k=1, 2, \ldots, \quad n=0, 1, \ldots,
\end{equation}
with initial conditions
\begin{equation}
  U_1^{n} = 0, \quad U_2^{n} = n, \quad  U_3^{n}=S_n,
  \qquad n=0, 1, \ldots.
  \label{am2}
\end{equation}
We have
\[
T_k^{(n)} = U_{3k+3}^{n}.
\]

\section{Numerical experiments}
\label{Sec:NumExp}

In this section, we will show how the convergence of some of some slowly
convergent example sequences can be accelerated by the new algorithm
\eqref{am1}--\eqref{am2}.

\noindent {\bf Example 1.}~We consider the linearly convergent sequence
\begin{eqnarray*}
  S_n = 2^{n} \sin \left(\frac{\pi}{2^n} \right),
\end{eqnarray*}
which converges to $S=\pi=3.14159~26535~89793\ldots$. The
corresponding transformation results are presented in Table 5.1.

\begin{table}
\caption{Numerical results of example 1.}
  \centering
  \begin{tabular}{|c|c|c|c|c|c|}\hline
    $n$ & $T_0^{(n)}$ & $T_1^{(n)}$ & $T_2^{(n)}$ & $T_3^{(n)}$ &
    $T_4^{(n)}$\\\hline
    1 & 2.00000 & 3.16790~51916 & 3.14158~12622 &
    3.14159~26537 & 3.14159~26536 \\
    2 & 2.82843 & 3.14304~69467 & 3.14159~24821 &
    3.14159~26536 & \\
    3 & 3.06146 & 3.14168~10168 & 3.14159~26509
    & & \\
    4 & 3.12145 & 3.14159~81382 & 3.14159~26535
    & & \\
    5 & 3.13655 & 3.14159~29958 & 3.14159~26536 & &
    \\\hline
  \end{tabular}
\end{table}

\noindent{\bf Example 2.}~We consider the sequence
\begin{eqnarray*}
  S_n = \sum_{k=1}^{n} \frac{(-1)^{k-1}}{k}
\end{eqnarray*}
of the partial sums of the alternating series, which converges to $S=\ln2=0.69314~71805\ldots$. The
corresponding transformation results are presented in Table 5.2.

\begin{table}
\caption{Numerical results of example 2.}
  \begin{tabular}{|c|c|c|c|c|c|c|}\hline
    $n$ & $T_0^{(n)}$ & $T_1^{(n)}$ & $T_2^{(n)}$ & $T_3^{(n)}$ &
    $T_4^{(n)}$ & $T_5^{(n)}$\\\hline
    1 & 1.00000 & 0.70588 & 0.69381 & 0.69318~36537 & 0.69314~92236 & 0.69314~72961\\
    2 & 0.50000 & 0.68817 & 0.69294 & 0.69313~74926 & 0.69314~66861 & 0.69314~71544\\
    3 & 0.83333 & 0.69557 & 0.69323 & 0.69315~03633 & 0.69314~73258 & 0.69314~71876\\
    4 & 0.58333 & 0.69178 & 0.69311 & 0.69314~59611 & 0.69314~71311 & \\
    5 & 0.78333 & 0.69399 & 0.69316 & 0.69314~77054 & 0.69314~71995 & \\
    \hline
  \end{tabular}
\end{table}

\noindent{\bf Example 3.}~We consider the logarithmically convergent
sequence
\begin{eqnarray*}
  S_n=\sum_{k=1}^{n}\frac{1}{k^2}
\end{eqnarray*}
of the partial sums of the Dirichlet series for $\zeta (2)$, which
converges to $S=\frac{\pi^2}{6}=1.64493~406684\ldots$. The corresponding
transformation results are presented in Table 5.3.

\begin{table}
\caption{Numerical results of example 3.}
  \begin{tabular}{|c|c|c|c|c|c|c|c|c|}\hline
    $n$ & $T_0^{(n)}$ & $T_1^{(n)}$ & $T_2^{(n)}$ & $T_3^{(n)}$ &
    $T_4^{(n)}$ & $T_5^{(n)}$ &$T_6^{(n)}$&$T_7^{(n)}$\\\hline
    1 & 1.00000 & 1.38462 & 1.49536 & 1.54487 & 1.57198 & 1.58872 & 1.59990 & 1.60782\\
    2 & 1.25000 & 1.45686 & 1.52776 & 1.56266 & 1.58298 & 1.59608 & 1.60511 & 1.61164\\
    3 & 1.36111 & 1.49794 & 1.54871 & 1.57512 & 1.59112 & 1.60175 & 1.60925 & 1.61486\\
    4 & 1.42361 & 1.52436 & 1.56334 & 1.58432 & 1.59738 & 1.60625 & 1.61261 & 1.61721 \\
    5 & 1.46361 & 1.54276 & 1.57412 & 1.59138 & 1.60234 & 1.60990 & 1.61542 & 1.62011\\\hline
  \end{tabular}
\end{table}

\noindent {\bf Note}

\noindent The results given above show that our new algorithm accelerates
convergence in the first two cases but not in the third case. As is well
known \cite{3}, Wynn's $\varepsilon$-algorithm, which corresponds to the
discrete potential KdV equation, can accelerate linear convergence and
the convergence of alternating series but fails to accelerate logarithmic
convergence. Apparently, the discrete Boussinesq equation and the
discrete potential KdV equation have similar properties from the
integrable systems' point of view. We therefore guess that the algorithm
based on \eqref{ag1} related to the lattice Boussinesq equation has the
similar acceleration properties as the $\varepsilon$-algorithm. The
numerical results given above are in agreement with our conjecture.

\section{Conclusions and discussions}
\label{Sec:ConclDisc}

In this article, we construct the molecule solution of \eqref{ag1}
related to the lattice Boussinesq equation by Hirota's bilinear
method. We show that this equation can be used as a numerical convergence
acceleration algorithm.  Our numerical experiments show that this
algorithm is effective for linearly convergent sequences and alternating
series but fails in the case of logarithmic convergence. Also we have
studied the confluent form of \eqref{ag1} in \cite{15b}.

It is known that the fully discrete potential KdV equation and the
lattice Boussinesq equation are the first two equations of the
Gel'fand--Dikii hierarchy. Now we know that both of the equations have
connections with convergence acceleration algorithms. Therefore, it is a
natural idea to investigate further whether higher order members of the
Gel'fand--Dikii hierarchy have relationships to other convergence
acceleration algorithms. We will consider these problems in the future.

\section*{Acknowledgments}

This work was partially supported by the National Natural Science
Foundation of China (Grant no. 11071241) and the knowledge innovation
program of LSEC and the Institute of Computational Mathematics, AMSS,
CAS. E.J.W. gratefully acknowledges the hospitality of the Institute of
Computational Mathematics of the Chinese Academy of Sciences where a part
of this work was done.

\end{document}